\documentclass[12pt]{article}
\usepackage{amsmath,amssymb,amsfonts,amscd}

\newtheorem{thm}{Theorem}[section]
\newtheorem{cor}[thm]{Corollary}
\newtheorem{lem}[thm]{Lemma}
\newtheorem{prop}[thm]{Proposition}
\newtheorem{df}[thm]{Definition}

\title {On 0-homology of categorical at zero semigroups}
\author {~B.~V.~\,Novikov, L.~Yu.~\,Polyakova \footnote{
Supported by N.~I.~Akhiezer foundation grant.}}

\begin{document}

\date{}
\maketitle

\begin{center}
\begin{minipage}[c]{13.5cm}
\footnotesize \noindent \textbf{Abstract.} The isomorphism of 0-homology groups of a
categorical at zero semigroup and homology groups of its 0-reflector is proved. Some
applications of 0-homology to Eilenberg--MacLane homology of semigroups are given.

\end{minipage}
\end{center}

\medskip

Homology of semigroups and monoids was defined in works by Eilenberg and MacLa\-ne, but
it was not developed well later on and turned out to be less investigated than cohomology
of semigroups. Nevertheless they find their application for different problems.

For instance, it is well-known~\cite{CartEilOrig} that if a group $G$ is a group of
fractions of its submonoid $M$ then $H_n(G,A)\cong H_n(M,A)$ for every $G$-module $A$. In
such a situation Dehornoy and Lafont \cite{DehLaf03} construct free resolutions for
monoids which allow, in particular, to compute the homology of braid groups.

In~\cite{Squier87} Squier showed that every monoid, possessing a finite complete
rewriting system, satisfied some homological condition. He answered negatively the
question on an existence of such a system for every finitely presented monoid with the
solvable word problem. Squier's approach was developed, for instance, in~\cite{Kob90}
where the method to construct a free resolution for monoids with a complete rewriting
system was described, which allowed in its turn to find homology of such monoids.

Homology of free partially commutative monoids arise in the articles by Husainov
(see~\cite{HusTkEng},~\cite{Hus04}) in connection with the construction of homology
groups of asynchronous transition systems.

If a semigroup $S$ contains the zero then its homology and
cohomology are trivial. In \cite{Novikov78Eng} (see also
\cite{Novikov2001}) so called 0-homology was built which, generally
speaking, is nontrivial for semigroups with zero. Furthermore, if
$S$ contains the zero then a semigroup $\bar{S}$, called 0-reflector
of $S$, can be constructed (see definition below) and its cohomology
groups are closely connected with 0-cohomology groups of $S$.
Moreover for categorical at zero semigroups these groups are
isomorphic in all dimensions. In particular 0-cohomology allows to
compute cohomology groups in some cases.

In view of the aforesaid in~\cite{0-homPart1Eng} 0-homology of semigroups is constructed
and it is shown that the properties of the first 0-homology groups are similar to those
of the first 0-cohomology groups. In this work we study 0-homology groups of greater
dimensions. The main result about an isomorphism of 0-homology groups of the categorical
at zero semigroup and  homology groups of its 0-reflector is contained in Section 2
(notice that its proof essentially differs from the proof of the similar statement for
0-cohomology). Section 3 is devoted to the defining relations of categorical at zero
semigroups and is auxiliary (however is of its own interest). It is used in Section 4 in
examples and applications of 0-homology to the computation of Eilenberg--MacLane homology
groups.

\section{Preliminaries and basic definitions}

All the modules under consideration are right modules. The notion
${\sf Sem}$ is used for the category of semigroups.

Consider ${\sf Sem}_0$ --- the category, which objects are
semigroups with zero elements, and morphisms are such mappings
$\varphi:S \to T$ that $\varphi(0)=0$, $\varphi^{-1}(0)=0$ and
$\varphi(xy)=\varphi(x)\varphi(y)$ if $xy\neq 0$ ({\it
0-homomorphisms}). The subcategory of the category  ${\sf Sem}_0$
consisting of semigroups with an adjoint zero is, obviously,
isomorphic to the category ${\sf Sem}$. Therefore we will consider
${\sf Sem}$ as the subcategory in ${\sf Sem}_0$.

Recall the definition of a reflective subcategory.

\begin{df}{\rm \cite{MacLaneCat}} \label{df_0-reflect} A
subcategory ${\sf D}$ of a category ${\sf C}$ is called {\it
reflective} if to each object $C \in {\sf C}$ such an object
$R_{\sf D}(C) \in {\sf D}$ (called {\sf D}-{\it reflector} of $C$)
and a morphism  $\varepsilon_{\sf D}(C): C\rightarrow R_{\sf
D}(C)$ are assigned that for every $D \in {\sf D}$ the diagram
\[
\begin{array}{ccc}
C&\mathop{\longrightarrow}\limits^{\varepsilon_{\sf D}(C)}&R_{\sf D}(C)\\
\downarrow & & \\
D& &
\end{array}\]
can be uniquely completed by a morphism from ${\rm Hom}_{{\sf
D}}(R_{\sf D}(C),D)$ up to commutative one.
\end{df}

In \cite{0-homPart1Eng} it was shown that the category  ${\sf Sem}$
is reflective in ${\sf Sem}_0$. For a semigroup $S$ with zero we
denote its ${\sf Sem}$-reflector by $\bar{S}$ and  call it a
0-reflector. For instance, the semigroup $S$ with an adjoint zero
has a trivial 0-reflector: $\bar{S}=S$.

The semigroup $\bar{S}$ admits other equivalent constructions. Let $S$ be given by
nonzero generators and defining relations:

\begin{equation}\label{kopred}
    S=\langle a_1,\dots, a_k \mid P_i=Q_i, 1\leq i\leq n \rangle.
\end{equation}
We say that a relation $P_i=Q_i$ is a {\it zero relation} if the value of the word $P_i$
in the semigroup $S$ equals 0.

\begin{prop}\label{prop_kopred} {\rm \cite{Novikov83Eng, 0-homPart1Eng}}.
If all the zero relations in  (\ref{kopred}) are thrown off then the semigroup obtained
is the 0-reflector of $S$.
\end{prop}

The following construction~\cite{Novikov78Eng, Novikov83Eng} is
convenient for direct work with the elements of the semigroup
$\bar{S}$.

Let $(S)$ denote the set of all sequences $(s_1,s_2,\dots ,s_n)$,
$n\geq 1$, for which the following conditions hold:
\begin{align*}
& s_i\in S\setminus 0 \mbox{ for every } 1\leq i\leq n; \\
& s_is_{i+1}=0 \mbox{ for every } 1\leq i\leq n-1.
\end{align*}

Define on the set $(S)$ such a binary relation $\nu$ that $(s_1,\dots ,s_m)\nu(t_1,\dots
,t_n)$ if and only if one of the following conditions holds:

1) $m=n$ and there exists $i$ ($1\leq i\leq m-1$) such that $s_i=t_iu,\ t_{i+1}=us_{i+1}$
for some $u\in S$ and $s_j=t_j$ if $j\ne i$, $j\ne i+1$;

2) $m=n+1$ and there exists $i$ ($2\leq i\leq m-1$) such that
$s_i=uv$, $t_{i-1}=s_{i-1}u$, $t_i=vs_{i+1}$ for some $u,v\in S$,
and $s_j=t_j$ if $1\leq j\leq i-2$, and $s_j=t_{j-1}$ if $i+2\leq
j\leq m$.

Let $\bar{\nu}$ be the least equivalence containing $\nu$ and
$\bar{S}$ be a quotient set $(S)/\bar{\nu}$.  Let $\langle s_1,
\dots , s_n\rangle $ denote the image of the element $(s_1,\dots
,s_n)\in (S)$ under factorization. Then $\bar{S}$ becomes a
semigroup which elements are multiplied by the following rule:
\[
\langle s_1,\dots ,s_m\rangle \langle t_1,\dots ,t_n\rangle =
\left\{
\begin{array}{ll}
    \langle s_1,\dots ,s_{m-1}, s_mt_1, t_2,\dots ,t_n \rangle  & \hbox{ if } s_mt_1\ne 0; \\
    \langle s_1,\dots ,s_{m-1}, s_m, t_1, t_2,\dots ,t_n \rangle  & \hbox{ if } s_mt_1=0.
\end{array}
\right.
\]

The following notions give us an Abelian category required for building 0-homology.

\begin{df} {\rm \cite{Novikov78Eng}}
A 0-module over a semigroup $S$ with zero is an Abelian group $A$ (in additive notation)
with a multiplication $A\times (S\setminus 0)\rightarrow A$ satisfying for every $s,t\in
S\setminus 0$, $a,b\in A$ the conditions

1) $(a+b)s=as+bs$,

2) if $st\ne 0$ then $(as)t=a(st)$.

A homomorphism from a 0-module $A$ to a 0-module $B$ (over $S$) is an Abelian groups
homomorphism $\varphi:A\rightarrow B$ such that $\varphi(as)=\varphi(a)s$ for every $s\in
S\setminus 0$, $a\in A$.
\end{df}

0-Modules over the semigroup $S$ form the category ${\sf C}_0(S)$ that is isomorphic to
the category ${\sf C}(\bar{S})$ of ordinary modules over $\bar{S}$ \cite{Novikov78Eng}.
The correspondence between objects of these categories is specified in such a way. If
$A\in {\sf C}_0(S)$ then $A$ becomes an $\bar{S}$-module by putting $a\langle s_1,\dots ,
s_n\rangle =(\dots(as_1)s_2\dots)s_n$ for $a\in A$. If $\bar{A}\in {\sf C}(\bar{S})$ then
$\bar{A}$ can be transformed into a 0-module over $S$ by putting $\bar{a}s=\bar{a}\langle
s\rangle $ for $\bar{a}\in \bar{A}$ and $s\in S\setminus 0 $.

Define now  0-homology groups for a semigroup $S$ with zero 0~\cite{0-homPart1Eng}. Let
$A$ be a 0-module over $S$. By $D_n$ we denote the subset of all $n$-tuples $[s_1,\dots,
s_n]$, where $s_j\in S$, $j=1,\dots , n$, such that $s_1s_2 \dots s_n\ne 0$. Let
$C_n^0(S, A)$ ($n\geq 1$) denote the set of all (finite) linear combinations of elements
from $D_n$ with coefficients in $A$. We write down such a linear combination as $\sum
a_{s_1,\dots ,s_n}[s_1,\dots ,s_n]$ and call it an $n$-dimensional $0$-chain. We put
$C_0^0(S,A)=A$.

The sets $C_n^0(S,A)$  ($n\geq 0$) are Abelian groups with respect to  addition. Define
boundary homomorphisms $\partial_n:C_n^0(S,A)\rightarrow C_{n-1}^0(S,A)$ on the
generators in an usual way:
\begin{align*}
\partial_n(a[s_1,\dots,
s_n]) = & as_1[s_2,\dots,s_n]+\sum_{i=1}^{n-1}(-1)^i
a[s_1,\dots,s_is_{i+1}, \dots ,s_n]+ \\
&(-1)^n a[s_1,\dots, s_{n-1}], \mbox{ if } n\geq 2;\\
\partial_1(a[s])=as-a.
\end{align*}

It is easy to see that $\partial_n$ is well defined  and is a boundary homomorphism:
$\partial_{n-1}\partial_n=0$.

\begin{df}
The group $H_n^0(S,A)={\rm Ker}\,\partial_n/ {\rm Im}\,
\partial_{n+1}$, $n\geq 1$ is called an $n$-th 0-homology group of
a semigroup $S$ with coefficients in a 0-module $A$.
\end{df}

In other words the 0-homology groups  $H_n^0(S,A)$ are defined as
the homology groups of the complex  $C_{\ast}^0$:
\[
\dots \stackrel{\partial_3}\rightarrow C_2^0(S,A)
\stackrel{\partial_2}\rightarrow C_1^0(S,A)
\stackrel{\partial_1}\rightarrow A.
\]

Along with the 0-homology of $S$ we consider homology groups $H_n(\bar{S},A)$ of the
0-reflector $\bar{S}$ with coefficients in the module $A$. According to one of the
definitions \cite{CartEilOrig} they are the homology groups of the complex $C_{\ast}$:
\[
\dots \stackrel{\delta_3}\rightarrow C_2(\bar{S},A)
\stackrel{\delta_2}\rightarrow C_1(\bar{S},A)
\stackrel{\delta_1}\rightarrow A.
\]
Here $C_j(\bar{S},A)$ are groups of chains, i.e. linear combinations of the form
$$
\sum a_{x_1,\dots,x_n}[x_1,\dots,x_n]
$$
where $a_{x_1,\dots,x_n}\in A$, $[x_1,\dots,x_n]$ are all the possible $n$-tuples of
$\bar{S}$ elements and only fi\-nitely many summands are nonzero. The boundary
homomorphisms $\delta_n$ are similar to the 0-chains homomorphisms $\partial_n$.

\textbf{Remark.} The zero homology group $H_0(\bar{S},A)$ equals
$A/{\rm Ker}\,\delta_1$, ${\rm Ker}\,\delta_1$ being generated by
all differences of the form $at-a$ where $a\in A$, $t\in \bar{S}$.
Therefore it is natural to define a zero 0-homology group as
\[ H_0^0(S,A)=A/{\rm Ker}\,\partial_1.\]
The group ${\rm Ker}\,\partial_1$ is a subgroup in $A$ generated by all differences of
the form $as-a$ where $a\in A$, $s\in S\setminus 0$. The equality $as-a=a\langle s\rangle
-a$ induces the embedding ${\rm Ker}\,\partial_1\hookrightarrow {\rm Ker}\,\delta_1$.
Since each generator  $a\langle s_1, \dots ,s_m \rangle-a$ of the group ${\rm
Ker}\,\delta_1$ can be represented as
\begin{align*}
& a\langle s_1, \dots ,s_m \rangle-a=(\dots (as_1)\dots )s_n-a = \\
& \Bigl((\dots (as_1)\dots s_{n-1})s_n-(\dots (as_1)\dots )s_{n-1}\Bigr)+ \\
& \Bigl( (\dots (as_1)\dots s_{n-2})s_{n-1}-(\dots (as_1)\dots
)s_{n-2} \Bigr) +   \dots + \\
& \Bigl((as_1)s_2-as_1\Bigr)+\Bigl(as_1-a\Bigr),
\end{align*}
this embedding is surjective. Hence, ${\rm Ker}\,\partial_1={\rm
Ker}\,\delta_1$ and
\[
H_0(\bar{S}, A)=H_0^0(S,A).
\]

Notice that if $S$ is a semigroup with the adjoint zero then
$H_n^0(S, A)\cong H_n(S\setminus 0,A)$.

If we consider $A$ as a 0-module over the semigroup $S$ and as an ordinary module over
$\bar{S}$ then the group homomorphism $\varepsilon_n:C_n^0(S,A)\rightarrow C_n(\bar{S},
A)$ ($n\geq 1$), defined as
\[ \varepsilon_n(a[s_1,\dots, s_n])=a[\langle
s_1\rangle ,\dots, \langle s_n\rangle ],
\]
arises in a natural way.

We put also $\varepsilon_0={\rm id}_A$. The homomorphisms family
$\varepsilon=\{\varepsilon_n\}_{n=0}^\infty$ can be represented as
a complex map:

\begin{picture}(250,60)(-100,-20)
\put(0,20){\dots}

\put(20,25){\vector(1,0){20}}  \put(45,20){$C_2^0(S,A)$}
\put(95,25){\vector(1,0){20}}  \put(120,20){$C_1^0(S,A)$}
\put(170,25){\vector(1,0){20}} \put(195,20){$A$}

\put(0,-20){\dots}

\put(20,-15){\vector(1,0){20}}  \put(45,-20){$C_2(\bar{S},A)$}
\put(95,-15){\vector(1,0){20}}  \put(120,-20){$C_1(\bar{S},A)$}
\put(170,-15){\vector(1,0){20}} \put(195,-20){$A$}

\put(65,14){\vector(0,-1){20}}  \put(140,14){\vector(0,-1){20}}
\put(200,14){\vector(0,-1){20}}

\scriptsize

\put(70, 2){$\varepsilon_2$} \put(145, 2){$\varepsilon_1$}
\put(205, 2){$\varepsilon_0$}

\put(24,29){$\partial_3$} \put(99,29){$\partial_2$}
\put(174,29){$\partial_1$}

\put(24,-11){$\delta_3$} \put(99,-11){$\delta_2$}
\put(174,-11){$\delta_1$}

\normalsize
\end{picture}

It is not difficult to check that for $i\geq 1$  the equalities
$\varepsilon_{i-1}\partial_i=\delta_{i}\varepsilon_i$ hold. Hence, the map
$\varepsilon=\{\varepsilon_n\}_{n=0}^\infty$ is a chain map.

Thus the homomorphisms $\varepsilon_n$ induce the homomorphisms
$\varepsilon_n^{\ast}:H_n^0(S,A)\rightarrow H_n(\bar{S}, A)$. For an arbitrary semigroup
$S$ with zero and a 0-module $A$ the following result was obtained in
\cite{0-homPart1Eng} (taking into account the remark, given above, about the isomorphism
of zero homology groups):

\begin{thm}\label{thm_1izom}
$\varepsilon_k^{\ast}$ is an isomorphism for $k\leq 1$ and an
epimorphism for $k=2$.
\end{thm}

\section{The main theorem}

In what follows $A$ is a fixed 0-module over a semigroup  $S$
unless specified otherwise.

\begin{df}{\rm \cite{CliffPrOrig2}}
A semigroup $S$ is called categorical at zero if $xyz=0$ implies
$xy=0$ or $yz=0$.
\end{df}

Our main result is contained in the following theorem:

\begin{thm} \label{thm_Kat0Izom}
If $S$ is categorical at zero then the map
$\varepsilon_n^{\ast}:H^0_n(S,A)\rightarrow H_n(\bar{S}, A)$ is an
isomorphism for all $n\geq 0$ and every 0-module $A$.
\end{thm}

In view of Theorem~\ref{thm_1izom} the statement has to be proved
only for $n\geq 2$. To prove that $\varepsilon_n^{\ast}$ is a
monomorphism we make use of the following Lemma:

\begin{lem}\label{lem_osn}{\rm \cite{0-homPart1Eng}} Let chain complexes
$M$, $N$ and a chain map $\alpha:M\rightarrow N$ be given:

\begin{picture}(250,65)(-100,-25)
\put(0,20){\dots} \put(200,20){\dots}

\put(20,25){\vector(1,0){20}}  \put(45,20){$M_{k+1}$}
\put(75,25){\vector(1,0){20}}  \put(100,20){$M_{k}$}
\put(120,25){\vector(1,0){20}} \put(145,20){$M_{k-1}$}
\put(175,25){\vector(1,0){20}}

\put(0,-20){\dots} \put(200,-20){\dots}

\put(20,-15){\vector(1,0){20}}  \put(45,-20){$N_{k+1}$}
\put(75,-15){\vector(1,0){20}}  \put(100,-20){$N_{k}$}
\put(120,-15){\vector(1,0){20}} \put(145,-20){$N_{k-1}$}
\put(175,-15){\vector(1,0){20}}

\put(53,14){\vector(0,-1){20}}  
\put(108,14){\vector(0,-1){20}} 
\put(153,14){\vector(0,-1){20}} 

\scriptsize

\put(75,30){$\partial_{k+1}$} \put(124,30){$\partial_k$}
\put(75,-10){$\delta_{k+1}$}  \put(124,-10){$\delta_k$}

\put(33, 2){$\alpha_{k+1}$}  

\put(92, 2){$\alpha_{k}$}    

\put(133, 2){$\alpha_{k-1}$} 

\normalsize
\end{picture}

If for some $k\geq 1$ there exist module homomorphisms
$\beta_j:N_j\rightarrow M_j$ ($j=k,k+1$) such that
\begin{align}
 & \beta_k \alpha_k -{\rm id}_{M_k}=0, \label{usl1}\\
 & \partial_{k+1}\beta_{k+1}=\beta_k\delta_{k+1}, \label{usl2}
\end{align}
then the induced homology groups homomorphism
$\alpha_k^{\ast}:H_k(M)\rightarrow H_k(N)$ is a mono\-mor\-phism.
\end{lem}

We put $M_k=C_k^0$, $N_k=C_k$, $\alpha_k=\varepsilon_k$ and
construct suitable homomorphisms $\beta_k$.

The following notations will be convenient: let $X_{l_i}^i$ denote
an element $\langle x_1^i, \dots , x_{l_i}^i\rangle \in \bar{S}$.
Besides if for an $n$-dimensional chain $a[X_{l_1}^1,\dots
,X_{l_n}^n]\in C_n(\bar{S}, A)$ the conditions $l_j\ne 1$,
$l_{j+1}=\dots =l_{n-1}=1$ hold, we put $x^j=x_{l_j}^j$,
$x^i=x^i_1$, $i=j+1,\dots, n$.

Define homomorphisms $\beta_n$ for $n\geq 2$ on the generators of groups $C_n(\bar{S},
A)$
\begin{align*}
& \beta_n(a[X_{l_1}^1,\dots ,X_{l_n}^n])= \left\{
\begin{array}{ll}
    aX_{l_1-1}^1[x^1,x^2, \dots ,x^n],  \hbox{ if } l_2=\dots=l_{n-1}=1 \\
    \qquad \qquad \qquad \ \ \ \ \ \ \ \ \ \ \ \hbox{ and } x^1x^2\dots x^n\ne 0; \\
    0,  \hbox{ otherwise.} \\
\end{array}
\right.
\end{align*}
and then extend them by linearity. In the proofs of the following two lemmas it is
sufficient to verify identities on the generators of corresponding groups. That is what
we will use.

\begin{lem}\label{lem_zetaeps}
For $n\geq 2$ the equality $\beta_n\varepsilon_n={\rm id}_{C_n^0}$
holds.
\end{lem}
\textbf{Proof.} For $n\geq 2$ and $a[s_1,\dots, s_n]\in C_n^0$ we
have:
\[
\beta_n\varepsilon_n(a[s_1,\dots, s_n])=\beta_n(a[\langle s_1\rangle
,\dots, \langle s_n\rangle ])=a[s_1,\dots, s_n],
\]
since $s_1s_2\dots s_n\ne 0$.  $\square$

\begin{lem} \label{lem_KommutKvKat0}
Let $S$ be a categorical at zero semigroup. Then
$\partial_n\beta_n=\beta_{n-1}\delta_n$ for all $n\geq 2$.
\end{lem}
\textbf{Proof.} For $n=2$ Lemma~\ref{lem_KommutKvKat0} is a
special case of Lemma 2.5 from \cite{0-homPart1Eng}. Let $n\geq
3$. Consider three possible cases for a generator
$a[X_{l_1}^1,\dots, X_{l_n}^n]\in C_n(\bar{S}, A)$.

\textbf{1.} Let $l_2=\dots=l_{n-1}=1$ and $x^1x^2\dots x^n\ne 0$.
Then
\begin{align*}
& \partial_n\beta_n \left(a[X_{l_1}^1,\langle x^2\rangle ,\dots,
\langle x^{n-1}\rangle , X_{l_n}^n] \right)= \partial_n \left( aX_{l_1-1}^1[x^1,\dots, x^n] \right)= \\
& aX_{l_1}^1[x^2,\dots, x^n]-aX_{l_1-1}^1[x^1x^2,\dots, x^n] + \\
& \sum_{j=2}^{n-1}(-1)^jaX_{l_1-1}^1[x^1,\dots,x^jx^{j+1},\dots,
x^n]+ (-1)^naX_{l_1-1}^1[x^1,\dots, x^{n-1}] = \\
& \beta_{n-1}\left ( aX_{l_1}^1[\langle x^2\rangle , \dots, \langle
x^{n-1}\rangle , X_{l_n}^n] - \right. \\
& a[\langle x_1^1,\dots,x_{l_1-1}^1,x^1x^2\rangle ,\langle
x^3\rangle ,\dots,
\langle x^{n-1}\rangle ,X_{l_n}^n]  + \\
& \sum_{j=2}^{n-1}(-1)^ja[X_{l_1}^1, \langle x^2\rangle , \dots,
\langle x^jx^{j+1}\rangle ,\dots,\langle x^{n-1}\rangle , X_{l_n}^n ]+ \\
&  \left. (-1)^na[X_{l_1}^1, \langle x^2\rangle , \dots,\langle
x^{n-1}\rangle ]\right) = \beta_{n-1}\delta_n \left(
a[X_{l_1}^1,\langle x^2\rangle ,\dots, \langle x^{n-1}\rangle ,
X_{l_n}^n]\right).
\end{align*}

\textbf{2.} Let $l_2=\dots=l_{n-1}=1$ and $x^1x^2\dots x^n=0$.
Then
$$
\partial_n\beta_n
\left( a[X_{l_1}^1,\langle x^2\rangle ,\dots, \langle
x^{n-1}\rangle , X_{l_n}^n] \right)=0
$$
and
\begin{align*}
& \beta_{n-1}\delta_n \left( a[X_{l_1}^1,\langle x^2\rangle
,\dots, \langle x^{n-1}\rangle , X_{l_n}^n] \right)=\\
& = \beta_{n-1}\left( aX_{l_1}^1[\langle x^2\rangle ,\dots,
\langle x^{n-1}\rangle , X_{l_n}^n]- a[X_{l_1}^1\langle x^2\rangle
,\langle x^3\rangle , \dots, \langle x^{n-1}\rangle ,
X_{l_n}^n]+ \right. \\
& + \sum_{j=2}^{n-2}(-1)^ja[X_{l_1}^1,\langle x^2\rangle
,\dots,\langle x^j\rangle \langle x^{j+1}\rangle ,\dots,
\langle x^{n-1}\rangle , X_{l_n}^n]+\\
& + \left. (-1)^{n-1} a[X_{l_1}^1,\langle x^2\rangle ,\dots,
\langle x^{n-1}\rangle X_{l_n}^n]+ (-1)^{n}a[X_{l_1}^1,\langle
x^2\rangle ,\dots, \langle x^{n-1}\rangle ] \right).
\end{align*}

In this expression $\beta_{n-1}$ vanishes on all the summands of
the intermediate sum. Notice that, if $x^1x^2=0$ and $x^2x^3\dots
x^n=0$, then $\beta_{n-1}$ vanishes on the other summands too. If
$x^1x^2=0$ but $x^2x^3\dots x^n\ne 0$ then $\beta_{n-1}$ vanishes
on each summand of the last pair and
\[
\beta_{n-1}\left( aX_{l_1}^1[\langle x^2\rangle ,\dots, \langle
x^{n-1}\rangle , X_{l_n}^n]- a[X_{l_1}^1\langle x^2\rangle
,\langle x^3\rangle , \dots, \langle x^{n-1}\rangle ,
X_{l_n}^n]\right) =0,
\]
since $X_{l_1}^1\langle x^2\rangle =\langle
x_1^1,\dots,x_{l_1-1}^1, x^1, x^2\rangle $.

If $x^1x^2\ne 0$ then $x^2x^3\dots x^n=0$ since $S$ is categorical at zero. Therefore
$\beta_{n-1}$ vanishes on each summand of the first pair. If at the same time
$x^2x^3\dots x^{n-1}\ne 0$ then categoricity at zero implies $x^{n-1}x^n=0$. Hence,
$\langle x^{n-1}\rangle X_{l_n}=\langle x^{n-1},x^{n}, x_2^n,\dots, x_{l_n}^n\rangle $
and
\[
\beta_{n-1}\left((-1)^{n-1} a[X_{l_1}^1,\langle x^2\rangle ,...,
\langle x^{n-1}\rangle X_{l_n}^n]+ (-1)^{n}a[X_{l_1}^1,\langle
x^2\rangle ,..., \langle x^{n-1}\rangle ] \right)=0.
\]
Thus in this case  $\beta_{n-1}\delta_n=0$.

\textbf{3.} Finally consider the case $l_m>1$ for some $2\leq m\leq n-1$. Again we have
$\partial_n\beta_n \left( a[X_{l_1}^1,\dots, X_{l_n}^n]\right )=0$ and
\begin{align*}
& \beta_{n-1}\delta_n \left( a[X_{l_1}^1,\dots, X_{l_n}^n] \right)=\\
& =\beta_{n-1}\left( aX_{l_1}^1[X_{l_2}^2,\dots, X_{l_n}^n]-
a[X_{l_1}^1X_{l_2}^2,\dots, X_{l_n}^n]+ \right.\\
& +
\sum_{j=2}^{n-2}(-1)^ja[X_{l_1}^1,\dots,X_{l_j}^jX_{l_{j+1}}^{j+1},\dots,
X_{l_n}^n]+ \\
& + \left. (-1)^{n-1}a[X_{l_1}^1,\dots,
X_{l_{n-1}}^{n-1}X_{l_n}^n]-a[X_{l_1}^1,\dots, X_{l_{n-1}}^{n-1}]
\right).
\end{align*}
In this expression $\beta_{n-1}$ vanishes on the summands of the
intermediate sum independently on $m$. If $2<m<n-1$ then
$\beta_{n-1}$ maps to 0 other summands as well. If $m=2$ then
$\beta_{n-1}$ vanishes on each summand of the last pair. Also
$\beta_{n-1}$ either equals 0 on the both summands of the first
pair or maps their sum to 0. Similarly for $m=n-1$.

The lemma is proved. $\square$

\begin{lem}\label{lem_epi}
Let $H_k^0(S,A)\cong H_k(\bar{S},A)$ for a semigroup $S$ with
zero, $k\geq 1$ and an arbitrary 0-module $A$. Then the map
$\varepsilon_{k+1}^{\ast}:H_{k+1}^0(S,A)\rightarrow
H_{k+1}(\bar{S},A)$ is an epimorphism.
\end{lem}
\textbf{Proof.} Consider $\bar{S}$-module $A$ as a quotient module
of a free $\bar{S}$-module $F$ by a submodule $B$. Thus we have a
commutative diagram with exact lines:

\noindent \small
\begin{picture}(300,60)(-40,-20)
\put(-15,20){\dots}

\put(0,25){\vector(1,0){15}} \put(20,20){$H_{k+1}^0(S,F)$}
\put(80,25){\vector(1,0){15}} \put(100,20){$H_{k+1}^0(S,A)$}
\put(160,25){\vector(1,0){15}} \put(180,20){$H_k^0(S,B)$}
\put(230,25){\vector(1,0){15}} \put(250,20){$H_k^0(S,F)$}
\put(300,25){\vector(1,0){15}} \put(320,20){$H_k^0(S,A)$}
\put(370,25){\vector(1,0){15}} \put(390,20){\dots}

\put(-15,-20){\dots}

\put(0,-15){\vector(1,0){15}} \put(20,-20){$H_{k+1}(\bar{S},F)$}
\put(80,-15){\vector(1,0){15}} \put(100,-20){$H_{k+1}(\bar{S},A)$}
\put(160,-15){\vector(1,0){15}} \put(180,-20){$H_k(\bar{S},B)$}
\put(230,-15){\vector(1,0){15}} \put(250,-20){$H_k(\bar{S},F)$}
\put(300,-15){\vector(1,0){15}} \put(320,-20){$H_k(\bar{S},A)$}
\put(370,-15){\vector(1,0){15}} \put(390,-20){\dots}

\put(45, 15){\vector(0,-1){22}} \put(125, 15){\vector(0,-1){22}}
\put(195, 15){\vector(0,-1){22}} \put(265, 15){\vector(0,-1){22}}
\put(335, 15){\vector(0,-1){22}}

\footnotesize \put(128,2){$\varepsilon_{k+1}^{\ast}$}
\put(338,2){$\varepsilon_k^{\ast}$}
\end{picture}
\normalsize

\medskip
Since $F$ is a free module $H_k(\bar{S}, F)=H_{k+1}(\bar{S}, F)=0$ and so $H_k(\bar{S},
B)\cong H_{k+1}(\bar{S}, A)$. Under the conditions $H_k^0(S, B)\cong H_k(\bar{S}, B)$ and
$H_k^0(S, F)\cong H_k(\bar{S}, F)=0$. It follows from here that the map $H_{k+1}^0(S,
A)\rightarrow H_k^0(S, B)$ is an epimorphism. Finally commutativity of the diagram
implies that the map $\varepsilon_{k+1}^{\ast}:H_{k+1}^0(S, A)\rightarrow
H_{k+1}(\bar{S}, A) $ is an epimorphism as well. $\square$

\noindent \textbf{The proof of Theorem~\ref{thm_Kat0Izom}.}
Lemmas~\ref{lem_zetaeps},~\ref{lem_KommutKvKat0} imply
conditions~(\ref{usl1}),~(\ref{usl2}) in Lemma~\ref{lem_osn} for
$n\geq 2$ being satisfied. Hence, the map $\varepsilon_n^{\ast}$
is a monomorphism. Applying successively Lemma~\ref{lem_epi} for
$k=2,3,\dots$ we obtain the required statement. $\square$

\section{Defining relations of categorical at zero se\-mi\-groups}

Denote by $S=\langle a_1, \dots , a_n \mid A_i=B_i, i=1, \dots , r \rangle$ a semigroup
with generators $a_j$ ($1\leq j\leq n$) and defining relations  $A_i=B_i$, $i=1, \dots ,
r$. Let $S$ be categorical at zero. If some defining relation of the semigroup $S$ is of
the form $A=0$ then, in view of categoricity at zero, it is a consequence of some
equality $a_ia_j=0$. Therefore in what follows we suppose that on the set $N=\{1,2,\dots
, n\}$ a relation $\Gamma$ is given such that $(i,j)\in \Gamma \Leftrightarrow a_ia_j=0$
and we write down a categorical at zero semigroup in the following form:
\begin{equation}\label{form}
    S=\langle a_1, \dots , a_n \mid a_ia_j=0 \, \mbox {for }(i,j)\in \Gamma; A_k=B_k, k=1,
\dots , m \rangle,
\end{equation}
where $A_k\ne 0$ and $B_k\ne 0$ for all $k\leq m$.

Introduce the notations: $\Gamma a_i=\{a_j\mid (j,i)\in \Gamma
\}$, $a_i\Gamma =\{a_j \mid (i,j)\in \Gamma \}$. Besides denote
the length of a word $A$ by $l(A)$; we suppose that in
(\ref{form}) $l(A_k)\geq l(B_k)$ and $l(A_k)\geq 1$ for all $k\leq
m$.
\begin{prop}\label{prop_KritKat}
Let a semigroup $S$ be given in the form (\ref{form}). Let
$A_k=p_kA'_kq_k$ and $B_k=r_kB'_ks_k$. The words $A'_k$ and $B'_k$
can be empty and if $l(B_k)=1$ we suppose that $B_k=r_k=s_k$ (here
$p_k, q_k, r_k, s_k\in \{a_1,\dots , a_n\}$). The semigroup $S$ is
categorical at zero if and only if $\Gamma p_k=\Gamma r_k$ and
$q_k \Gamma =s_k \Gamma$ for all $k\leq m$.
\end{prop}
\textbf{Proof.} Let $S$ be categorical at zero and, for instance, $a\in \Gamma p_k$. Then
$A_k=B_k$ implies $ar_kB'_ks_k=0$. Since $B'_k\ne 0$, in view of categoricity, $ar_k=0$.
Hence, $a\in \Gamma r_k$.

Next verify the converse statement. Let $\Gamma p_k= \Gamma r_k$, $q_k\Gamma =s_k \Gamma$
and $XYZ=0$, $XY\ne 0$, $YZ\ne 0$. If the word  $XYZ$ contains a product $a_ia_j=0$ then
either $XY=0$ or $YZ=0$ which is impossible. Therefore the transformation of the word
$XYZ$ into zero is realized only by equalities $A_k=B_k$. However, according to the
condition, the products $a_ia_j=0$ cannot appear during such a transformation. Hence,
contrary to the assumption either $XY=0$ or $YZ=0$. $\square$

Consider now a connection between defining relations of semigroup $S$ and those of its
0-reflector $\bar{S}$.

\begin{prop}\label{prop_OpredSootnKat}
Let a categorical at zero semigroup $S$ be given by  defining
relations  (\ref{form}). Then
\[
\bar{S} =\langle \langle a_1\rangle, \dots , \langle a_n\rangle
\mid A_k=B_k, k=1, \dots , m \rangle,
\]
where the words  $A_k, B_k$ are considered in the alphabet
$\langle a_1\rangle, \dots , \langle a_n\rangle$. Conversely  if
$\bar{S}$ is given by relations $A_k=B_k$ ($k=1,\dots , m$) then
there exists a subset  $\Gamma \subseteq N$ such that the
semigroup $S$ can be given in the form (\ref{form}).
\end{prop}
\textbf{Proof.} The first part of the proposition follows immediately from
Proposition~\ref{prop_kopred}.

Let now the semigroup $\bar{S}$ be defined by the relations
 $A_k=B_k$ ($k=1,\dots , m$) and $C=D$ be an equality in
 $S$. If $C\not\equiv 0$, $D\not\equiv 0$ then this equality holds
in $\bar{S}$ as well. Hence it is a consequence of the relations $A_k=B_k$. If, for
instance,  $C\not\equiv 0$, $D\equiv 0$ and $C\equiv a_{i_1}\dots a_{i_r}$ then
categoricity at zero implies $a_{i_k}a_{i_{k+1}}=0$ for some $k$, i.e. the equality $C=0$
follows from $a_ia_j=0$, $(i,j)\in \Gamma$. The second part of the statement is proved.
$\square$

\section{Some applications}

The results of the previous section can be used to establish
connections between ordinary homology groups and 0-homology ones.
The following assertion is a simple example:
\begin{prop} \label{cor_NulOprSootn}
Let all the defining relations of a semigroup $S$ be of the form
$a_ia_j=0$. Then $H_n^0(S,A)=0$ for all $n>1$ and every 0-module
$A$ over $S$.
\end{prop}
\textbf{Proof.} According to Proposition~\ref{prop_KritKat} $S$ is
categorical at zero. Proposition~\ref{prop_OpredSootnKat} implies
that $\bar{S}$ is a free semigroup. Hence, $H_n(\bar{S}, A)=0$ for
$n>1$ (see, for example,~\cite{CartEilOrig}). Now the statement
follows from Theorem~\ref{thm_Kat0Izom}. $\square$

Usually a semigroup with zero is simpler than its 0-reflector.
Therefore for computation of homology groups of a given semigroup
$T$ the following technique can be used: find a categorical at zero
semigroup $S$ such that its 0-reflector $\bar{S}$ is isomorphic to
$T$, calculate $H_n^0(S,\_)$ and use Theorem \ref{thm_Kat0Izom}.

Let a semigroup $T$ be given in the form
\begin{equation}\label{form2}
    T=\langle a_1, \dots , a_n \mid A_k=B_k, k=1,\dots , m \rangle
\end{equation}

Introduce the notation:  ${\cal P}=\{A_k=B_k\mid 1\le k \le m\}$.
Let $I({\cal P})$ denote the set of the elements $x\in T$ such
that $A_k \not\in TxT$ for all $k$. This set is an ideal in $T$ if
it is not empty.

The following proposition proved in~\cite{0-homPart1Eng} is in
some sense the converse to Proposition~\ref{prop_kopred}. It will
be helpful for us in examples.
\begin{prop}\label{prop_ideal}
Let a semigroup $T$ be given in the form~(\ref{form2}) and
$I({\cal P})\ne \emptyset$. If  $a_j \not\in I({\cal P})$ for all
$1\leq j\le m$ then $T$  is a  0-reflector of the quotient
semigroup $S=T/I({\cal P})$.
\end{prop}

\textbf{Example 1.} Consider the semigroup $T=\langle a, b, c, d \mid ab=cd \rangle$.
Then $T\setminus I({\cal P})=\{a, b, c, d, x=ab=cd\}$ and by Proposition~\ref{prop_ideal}
$T=\bar{S}$ where $S$ consists of the elements $0, a, b, c, d, x$, all the products being
equal zero except $ab=cd=x$. Since $S^3=0$ we have $H_2^0(S,A)={\rm Ker}\,\partial_2$ for
each 0-module $A$ over $S$. An arbitrary 2-dimensional 0-cycle is of the form
$f=\alpha[a,b]+\beta[c,d]$. The equality $\partial_2 f=0$ implies $\alpha=\beta=0$ and so
$H_2^0(S,A)=0$. Hence $H_2(T,A)=0$ by Theorem~\ref{thm_1izom}. This implies $H_n(T,A)=0$
for every $n\ge 2$ and every $T$-module $A$.

\textbf{Example 2.} Let $T=\langle a, b, c \mid ab=ac \rangle$. In
this example $T\setminus I({\cal P})=\{a, b, c, ab\}$. Then by
Proposition~\ref{prop_ideal} $T=\bar{S}$ where $S$ consists of the
elements $0, a, b, c, ab$. Again $S^3=0$ and we have
$H_2^0(S,A)={\rm Ker}\,\partial_2$. Let $f=\alpha[a,b]+\beta[a,c]$
be a 2-dimensional 0-cycle. Then the equality $\partial_2 f=0$
implies
 $\alpha a=0$ and $\beta=-\alpha$. Thus the group
$H_2^0(S,A)$ is isomorphic to the subgroup $A_a$ of the 0-module
$A$ consisting of the elements $\alpha$ such that $\alpha a=0$. It
is not difficult to verify that $S$ is a  categorical at zero
semigroup. So by Theorem~\ref{thm_Kat0Izom} $H_2(T,A)\cong A_a$.

Let a semigroup  $T$ be given in the form~(\ref{form2}). Assign to
the defining relations system ${\cal P}$ the graph $\Delta$, which
vertices set $\{1,2,\dots , n\}$ and the edges are the pairs
$(i,j)$ such that the product $a_ia_j$ is contained in some of the
words $A_k, B_k$ ($k\leq m$).

We call a vertex $a$ of the graph $\Delta$ \emph{an entrance}
(\emph{an exit}) if $(b,a)\not\in \Delta$ (respectively
$(a,b)\not\in \Delta$) for every vertex $b$.
\begin{thm}\label{thm_graph}
Let a semigroup $T$, given in the form~(\ref{form2}), satisfy the
following condition: for all the words $A_k, B_k$ ($k\leq m$)
their first letters are entrances and the last letters are exits
in the graph $\Delta$.

Then the semigroup
\[
S=\langle x_1,\dots , x_n \mid x_ix_j=0 \Leftrightarrow
(i,j)\not\in \Delta ; A_k=B_k, k=1,\dots , m\rangle,
\]
where the words $A_k, B_k$ are written down in the alphabet
 $\{x_1,\dots , x_n\}$, is categorical at zero and
$T$ is a 0-reflector of $S$.
\end{thm}
\textbf{Proof.} Let $A_k=p_kA'_kq_k$, $B_k=r_kB'_ks_k$ where $p_k, q_k, r_k, s_k\in
\{x_1,\dots, x_n\}$. Since $p_k$ and $r_k$ are entrances $\Gamma p_k=\Gamma r_k=\{x_1,
\dots , x_n\}$ and similarly $q_k\Gamma=s_k\Gamma$. Because of the same reason none of
the words $A_k$, $B_k$, ($k\leq m$) contains the other. Therefore in the semigroup none
of the defining relations $A_k=B_k$ is a zero relation. According to
Proposition~\ref{prop_KritKat} $S$ is categorical at zero.
Proposition~\ref{prop_OpredSootnKat} implies that its 0-reflector is isomorphic to $T$.
$\square$

\begin{cor}\label{cor_graph} Let a semigroup $T$ is under the conditions of
the previous theorem and the graph $\Delta$ does not contain any circuit. Then
$H^l(T,A)=0$ for all $l> l_0+1$ where $l_0$ is the length of the longest path in
$\Delta$.
\end{cor}
\textbf{Proof.} Consider the semigroup $S$ from Theorem~\ref{thm_graph}. In consequence
of the absence of circuits every word in $S$ of the length greater than $l_0+1$ equals
zero and so $S^{l_0+2}=0$. Therefore $C_0^l(S,A)=0$ as soon as $l>l_0+1$. Hence,
$H_l(T,A)\cong H_l^0(S,A)=0$. $\square$

\textbf{Example 3.} Consider the Adyan semigroup
\cite{Adyan66Eng}:
\[
T=\langle a, b, c, d, e \mid ab=cd, aeb=ced \rangle.
\]
The graph $\Delta$ for it is of the form:

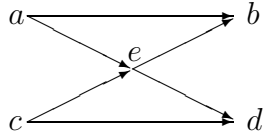
\begin{center}
\begin{picture}(120,60)(-60,-30)
\put(0,0){\vector(2,1){38}} \put(0,0){\vector(2,-1){38}}
\put(-40,20){\vector(2,-1){39}} \put(-40,20){\vector(1,0){78}}
\put(-40,-20){\vector(2,1){39}} \put(-40,-20){\vector(1,0){78}}

\put(-47, 17){$a$}\put(-47, -23){$c$}

\put(43, 17){$b$} \put(43, -23){$d$}

\put(-2,3){$e$}
\end{picture}
\end{center}

Therefore all homology groups of dimension 4 and greater are trivial.

In conclusion we consider free products of semigroups.

In~\cite{0-homPart1Eng} the description of the first homology group of free product of
two semigroups was obtained. If $S$ and $T$ are semigroups without zero and $A$ is a
$S\ast T$-module we denote by $A(S-1)$ (respectively by  $A(T-1)$) a subgroup in the
module $A$ generated by the elements of the form $as-a$ where $a\in A, s\in S$
(respectively $at-a$ where $a\in A, t\in T$) and put $A_1=A(S-1)\cap A(T-1)$. Then the
following proposition holds:
\begin{prop}
$H_1(S\ast T, A)$ is an extension of $H_1(S,A)\oplus H_1(T,A)$ by
$A_1$.
\end{prop}

In particular, if $A$ is a trivial $S\ast T$-module then $H_1(S\ast T, A)\cong
H_1(S,A)\oplus H_1(T,A) $.

Now proceed to the homology groups of greater dimensions. Recall
 \cite{CliffPrOrig2} that a semigroup $U$ is called a 0-direct union of semigroups
  $\{S_\lambda\}_{\lambda\in \Lambda}$ if
$U=\bigcup_{\lambda\in \Lambda} S_{\lambda}$, $S_\lambda\cap S_\mu=0$ and $S_\lambda
S_\mu=0$ for all $\lambda\ne \mu$.

\begin{lem}\label{lem_0-prObyed} Let a semigroup $S$ be a 0-direct
union of semigroups $S_{\lambda}$, $(\lambda\in \Lambda)$. Then
\[
H_n^0(S,A)\cong \bigoplus_{\lambda \in \Lambda}H^0_n(S_{\lambda},
A),
\]
where $A$ is a 0-module over $S$ (and so over every $S_\lambda$ as well) and $n>1$.
\end{lem}
\textbf{Proof.} Let $c=\sum a_{s_1\dots s_n}\otimes(s_1,\dots, s_n)$ be an
$n$-dimensional $0$-cycle from $C_n^0(S, A)$. Then for every summand $a_{s_1\dots
s_n}\otimes (s_1,\dots, s_n )$ all elements $s_j$ belong to the same semigroup
$S_{\lambda}$, otherwise $c$ is not defined. Therefore every cycle can be given as the
sum $c=\Sigma_{\lambda\in \Lambda}c_{\lambda}$ where $c_{\lambda}$ is a cycle belonging
to $C_n^0(S_{\lambda}, A)$. At the same time  $c=\partial c'$  for some $c'\in
C_{n+1}^0(S,A)$ if and only if $c_{\lambda}=\partial c'_{\lambda}$ for all
$\lambda\in\Lambda$, what implies the statement. $\square$

\begin{prop}\label{prop_SvProizv}
Let $S=\prod_{\lambda\in\Lambda}^{\ast}S_{\lambda}$ be a free
product of semigroups $S_\lambda$. Then
\[
H_n(S,A)\cong\bigoplus_{\lambda\in \Lambda}H_n(S_{\lambda}, A)
\]
for every $S$-module $A$  and each $n>1$.
\end{prop}
\textbf{Proof.} Similarly to the proof of Theorem~5 in \cite{Novikov78Eng} consider $T$
--- the 0-direct union of the semigroups  $T_{\lambda}=S_{\lambda}\cup 0$ with extra zeroes.
Then the semigroup $T$ is categorical at zero and $\overline{T}\simeq S$. In view of
Theorem~\ref{thm_Kat0Izom} and the previous lemma we obtain the required statement.
$\square$

\end{document}